\documentclass{amsart}

\usepackage{amsthm, amsfonts, amssymb, amsmath, latexsym, enumerate, times}
\usepackage[latin1]{inputenc}
\usepackage{mathrsfs}
\usepackage{mathtools}
\usepackage{letltxmacro}
\LetLtxMacro\orgvdots\vdots
\LetLtxMacro\orgddots\ddots

\makeatletter
\DeclareRobustCommand\vdots{%
	\mathpalette\@vdots{}%
}
\newcommand*{\@vdots}[2]{%
	\sbox0{$#1\cdotp\cdotp\cdotp\m@th$}%
	\sbox2{$#1.\m@th$}%
	\vbox{%
		\dimen@=\wd0 %
		\advance\dimen@ -3\ht2 %
		\kern.5\dimen@
		\dimen@=\wd2 %
		\advance\dimen@ -\ht2 %
		\dimen2=\wd0 %
		\advance\dimen2 -\dimen@
		\vbox to \dimen2{%
			\offinterlineskip
			\copy2 \vfill\copy2 \vfill\copy2 %
		}%
	}%
}
\DeclareRobustCommand\ddots{%
	\mathinner{%
		\mathpalette\@ddots{}%
		\mkern\thinmuskip
	}%
}
\newcommand*{\@ddots}[2]{%
	\sbox0{$#1\cdotp\cdotp\cdotp\m@th$}%
	\sbox2{$#1.\m@th$}%
	\vbox{%
		\dimen@=\wd0 %
		\advance\dimen@ -3\ht2 %
		\kern.5\dimen@
		\dimen@=\wd2 %
		\advance\dimen@ -\ht2 %
		\dimen2=\wd0 %
		\advance\dimen2 -\dimen@
		\vbox to \dimen2{%
			\offinterlineskip
			\hbox{$#1\mathpunct{.}\m@th$}%
			\vfill
			\hbox{$#1\mathpunct{\kern\wd2}\mathpunct{.}\m@th$}%
			\vfill
			\hbox{$#1\mathpunct{\kern\wd2}\mathpunct{\kern\wd2}\mathpunct{.}\m@th$}%
		}%
	}%
}
\makeatother

\usepackage{array}
\usepackage{comment}
\usepackage{paralist}
\usepackage{xcolor}

\newtheorem{theorem}{Theorem}
\newtheorem{lemma}[theorem]{Lemma}

\theoremstyle{definition}

\newcommand{\bbP}{{\mathbb P}}

\newcommand{\bbC}{{\mathbb C}}

\def\leq{\leqslant}
\def\le{\leqslant}
\def\ge{\geqslant}

\begin{document}
 
\title[On weighted Bounded Negativity]{A remark on weighted Bounded Negativity \\ for blow-ups of the projective plane}
 
\author{Ciro Ciliberto}
\address{Dipartimento di Matematica, Universit\`a di Roma Tor Vergata, Via O. Raimondo 00173 Roma, Italia}
\email{cilibert@axp.mat.uniroma2.it}

\author{Claudio Fontanari}
\address{Dipartimento di Matematica, Universit\`a degli Studi di Trento, Via Sommarive 14, 38123 Povo, Trento}
\email{claudio.fontanari@unitn.it}
 
\subjclass{Primary 14C17; Secondary 14J26}
 

\centerline{}
\begin{abstract} 
Motivated by the \emph{Weighted Bounded Negativity Conjecture} (see \cite{BB+}, Conjecture 3.7.1), 
we prove that all but finitely many reduced and irreducible curves $C$ on
the blow-up of $\mathbb{P}^2$ at $n$ points satisfy the inequality $C^2 \ge \min \{-\frac{1}{12} n (C.L +27), -2 \}$, 
where $L$ is the pull-back of a line. This partially improves on some result in \cite{LP}. 
\end{abstract}

\maketitle 

\section{Introduction} 

The celebrated \emph{Bounded Negativity Conjecture}, going back (at least) to Federigo Enriques (see for instance \cite{BH+}, Conjecture 1.1 and the historical remarks following its statement), predicts that on every smooth surface $S$ the self-intersection $C^2$ of any reduced and irreducible curve $C$ on $S$ is bounded below by a constant depending only on $S$. The extreme difficulty of such a conjecture, which is still widely open even for the blow-up of the projective plane 
$\bbP^2$ at $n \ge 10$ general points, motivated the formulation of a few weaker versions, among which the so-called \emph{Weighted Bounded Negativity Conjecture} 
(see \cite{BB+}, Conjecture 3.7.1): here the bound is allowed to depend on the degree of $C$ with respect to any nef and big divisor on $S$. The recent paper 
\cite{LP} collects several partial results towards this direction. In particular, in \cite{LP}, Theorem 3.1, by using Orevkov-Sakai-Zaidenberg's inequality, it is proven that if $Y$ is the blow-up of $\bbP^2$ at $n$ distinct points and $C$ is a reduced and irreducible curve on $Y$, then $C^2 \ge -2n C.L$, where $L$ is the pull-back of a line.\footnote{We guess that the above statement implicitly assumes $C.L > 0$, in order to exclude the exceptional divisors.} 
In the same paper (see p. 370), as a consequence of the Pl\"ucker-Teissier formula, the previous bound is improved to $C^2 \ge -n C.L$. 

Here, by applying only elementary tools in algebraic surface theory, we obtain the following partial improvement: 

\begin{theorem}\label{main}
Let $Y$ be the blow-up of $\bbP^2$ at $n$ points and let $L$ be the pull-back of a line in $\bbP^2$.
Then all but finitely many reduced and irreducible curves $C$ on $Y$ satisfy the inequality:
$$
C^2 \ge \min \{-\frac{1}{6} n (C.L +3), -2 \}.
$$
\end{theorem}

Next, as a consequence of Miyaoka-Yau-Sakai's inequality, we sharpen the previous bound as follows: 

\begin{theorem}\label{miyaoka}
Let $Y$ be the blow-up of $\bbP^2$ at $n$ points and let $L$ be the pull-back of a line in $\bbP^2$.
Then all but finitely many reduced and irreducible curves $C$ on $Y$ satisfy the inequality:
$$
C^2 \ge \min \{-\frac{1}{12} n (C.L +27), -2 \}.
$$
\end{theorem}

We work over the complex field $\bbC$.\medskip

{\bf Acknowledgements:} The authors are members of GNSAGA of the Istituto Nazionale di Alta Matematica ``F. Severi". This research project was partially supported by PRIN 2017 ``Moduli Theory and Birational Classification''. 

\section{The proofs}

\begin{lemma}\label{self-intersection}
Let $X$ be a smooth projective surface and let $C$ be a reduced and irreducible curve on $X$. Then for every integer $m \ne 1$ 
we have 
\begin{eqnarray*}
 C^2 = \frac{1}{m-1} \chi(\mathcal{O}_X) + \frac{1}{2} m K_X^2 + 2 p_a(C) + \frac{1}{m-1} p_a (C) - 2 - \frac{1}{m-1} \\
 - \frac{1}{m-1} h^0(mK_X+C) + \frac{1}{m-1} h^1(mK_X+C) - \frac{1}{m-1}h^0(-(m-1)K_X - C).
\end{eqnarray*}
\end{lemma}

\proof Just apply Riemann-Roch theorem to $mK_X+C$, Serre duality to $h^2(mK_X+C)$ and the adjunction formula to $C$. 
\qed

\begin{lemma}\label{anticanonical-nonvanishing}
Notation as in Lemma \ref{self-intersection}. If $h^0(-mK_X) \ne 0$ for some $m \ge 1$, then all but finitely many reduced and irreducible curves $C$ on $X$ satisfy the inequality $C^2 \ge -2$. 
\end{lemma}

\proof Let $E$ be an effective divisor in $\vert -mK_X \vert$. If $C$ is not one of the finitely many curves in the support of $E$ then 
$-K_X.C \ge 0$, hence $C^2 = -K_X.C + 2p_a(C) - 2 \ge -2$. 
\qed

\begin{lemma}\label{adjunction-extinction}
Notation as in Theorem \ref{main}. Let $m_0$ be the integer such that 
$$
\frac{C.L + 1}{3} \le m_0 \le \frac{C.L + 3}{3}.
$$
Then $h^0(m_0 K_Y + C) = 0$. 
\end{lemma}

\proof We have $(m_0 K_Y + C).L = -3 m_0 + C.L \leq -1 <0$ and $L^2 = 1 \ge 0$, hence the divisor $m_0 K_Y + C$ cannot be effective. 
\qed

\medskip

\noindent \textit{Proof of Theorem \ref{main}.}
We apply Lemma \ref{self-intersection} to $X = Y$, in particular we have $\chi(\mathcal{O}_Y)= 1$ and $K_Y^2 = 9-n$. 
By Lemma \ref{anticanonical-nonvanishing} we may assume $h^0(-(m-1)K_Y - C) = h^0(-(m-1)K_Y) = 0$. 
By setting $m = m_0$ as in Lemma \ref{adjunction-extinction} we obtain:
$$
C^2 \ge - \frac{1}{2} m_0 n \ge  -\frac{1}{6} n ( C.L +3).
$$
\qed

\noindent \textit{Proof of Theorem \ref{miyaoka}.}
If $h^0(2 K_Y + C) = 0$, then by arguing as in the proof of Theorem \ref{main} with $m=2$ we obtain $C^2 \ge \min \{-n, -2 \}$. 
Assume now $h^0(2 K_Y + C) \neq 0$, so that in particular $h^0(2 (K_Y + C)) \neq 0$. By Lemma \ref{anticanonical-nonvanishing} we may 
also assume $h^0(-K_Y) = 0$. Hence we are in the position to apply \cite{H}, Corollary 1.8, and deduce 
$$
C^2 \ge K^2_Y - 3 c_2(Y) + 2 - 2 p_a(C) = 9-n - 3(3+n) + 2 - 2 p_a(C) = -4n + 2 - 2 p_a(C).   
$$
On the other hand, by arguing as in the proof of Theorem \ref{main}, we obtain 
$$
C^2 \ge -\frac{1}{6} n ( C.L +3) + 2 p_a(C) - 2.  
$$
It follows that 
$$
2 C^2 \ge -4n + -\frac{1}{6} n ( C.L +3) 
$$
and we conclude 
$$
C^2 \ge -\frac{1}{12} n ( C.L +27).
$$
\qed

\end{document}